\documentclass[12pt]{article}
\usepackage{amssymb}
\usepackage{latexsym}
\def\part#1{\frac{\partial\phantom{q}}{\partial#1}}

\newenvironment{rmk}{\begin{trivlist}\item[]{\bf Remark:} }
{\end{trivlist}}
\newenvironment{rmks}{\begin{trivlist}\item[]{\bf Remarks:} }
{\end{trivlist}}
\newenvironment{ex}{\begin{trivlist}\item[]{\bf Example:} }
{\end{trivlist}}
\newenvironment{prf}{\begin{trivlist}\item[]{\bf Proof:} }
{\hfill $\Box$ \end{trivlist}}

\newtheorem{thm}{Theorem}

\newtheorem{prp}[thm]{Proposition}

\newcommand{\lie}[1]{\mathfrak{#1}}
\def\End{\mathop{\rm End}\nolimits}
\def\Hom{\mathop{\rm Hom}\nolimits}

\def\ker{\mathop{\rm ker}\nolimits}

\def\deg{\mathop{\rm deg}\nolimits}

\def\tr{\mathop{\rm tr}\nolimits}

\def\Res{\mathop{\rm Res}\nolimits}
\def\KK{{\mathrm K}}

\def\Nm{\mathop{\rm Nm}\nolimits}

\newcommand{\R}{\mathbf{R}}
\newcommand{\C}{\mathbf{C}}

\newcommand{\Z}{\mathbf{Z}}

\begin{document}
\title{Langlands duality and $G_2$ spectral curves}
 \author{Nigel Hitchin\\[5pt]}
\maketitle

\begin{abstract} We first demonstrate how duality for the fibres of the  so-called Hitchin fibration works for the Langlands dual groups $Sp(2m)$ and $SO(2m+1)$. We then show that duality for $G_2$ is implemented by an involution on the base space which takes one fibre to its dual. A formula for the natural cubic form  is given and shown to be invariant under the involution.
\end{abstract}

\section{Introduction}
A recent paper \cite{KW} of Kapustin and Witten described the geometric Langlands programme in terms of ${\mathcal N}=4$ super Yang-Mills theory. Within this picture a fundamental role is played by  Langlands  duality, originating in the duality between electric and magnetic charges described  by Montonen and Olive many years ago. This is the duality between  the root lattices of  Lie groups $G$ and   $^L\negthinspace
G$. 

These physical aspects have mathematical interpretations when applied, as in  Kapustin and Witten,  to a particular gauge-theoretic moduli space introduced by the author \cite{H1}, \cite{H2}. This is the moduli space of Higgs bundles with structure group $G$ on a Riemann surface $\Sigma$. A distinctive feature of this space is its interpretation as an algebraically completely integrable Hamiltonian system -- it is a holomorphic symplectic manifold with a proper map  to  a vector space, such that the generic fibre  is a complex Lagrangian torus which  is  an abelian variety. The duality then manifests itself in the statement that the dual of the abelian variety  for the group $G$ should be the abelian variety for the Langlands dual group $^L\negthinspace
G$. Hausel and Thaddeus in \cite{HT} have observed this fact in many cases, and obtained global results on the topology of the moduli spaces which reflect what is expected of mirror symmetry. A general proof of the result has been given by Donagi and Pantev \cite{DP}, and is implicitly to be found in \cite{F}.

Our purpose in this paper  is to describe in a concrete fashion this result for the special case of $G_2$. As remarked in the physics paper \cite{AKS} ``.... S-duality in the case of $F_4$ and $G_2$ acts nontrivially on the moduli space of the gauge theory..." What this means for us is that, although the Langlands dual of $G_2$ is again $G_2$, the dual of the abelian variety over one point in the base of the integrable system is the abelian variety over a {\it different} point. The base space for $G_2$ consists of pairs $(f,q)$ of differentials on $\Sigma$ of degrees $2$ and $6$ respectively, and it is the involution 
$$(f,q)\mapsto (f,\frac{1}{54}f^3-q)$$
which takes a fibre to its dual. 
We show that this involution preserves the natural  cubic form on the base. 

To analyse the duality,  we have to identify concretely the abelian variety for a $G_2$ Higgs bundle. Rather than appeal to the Lie-theoretical approach of \cite{DP}, or the earlier work of \cite{Katz}, we choose to describe $G_2$ as the identity component of the stabilizer of a three-form in seven dimensions, drawing on the papers \cite{H3},\cite{H4} dealing with $G_2$ in a different context.  The abelian variety is then described as  the intersection of two Prym varieties, and we can describe the duality in terms of the geometry of two spectral curves with a common double covering.

Since we use the 7-dimensional representation of $G_2$ in this approach, the Langlands duality between $Sp(6)$ and $SO(7)$ plays a role, and we take the opportunity here to describe the general case of duality for the  groups $Sp(2m)$ and $SO(2m+1)$. In the first half of the paper we provide a detailed description, replacing the sketchy version in the author's original  paper \cite{H2}, and at the same time correct an oversight, pointed out by Michael Thaddeus, whose resolution explains the duality. 
\vskip .15cm
The author wishes to thank M.Thaddeus, A.Kapustin and E.Witten for useful conversations, and the Clay Mathematics Institute and KITP for support during the preliminary stages of this work.

\section{The general linear group \label{standard}}
\subsection{Spectral curves}
A Higgs bundle for $G=GL(n)$ on a compact  Riemann surface $\Sigma$ of genus $g>1$ is a holomorphic vector bundle $E$ of rank $n$, together with a Higgs field $\Phi$, which is a holomorphic section of  $\End E\otimes K$ satisfying the stability condition that a $\Phi$-invariant subbundle has slope less than that of $E$. The moduli space ${\mathcal M}$ of such pairs $(E,\Phi)$ is the moduli space of Higgs bundles and is a complex orbifold of dimension $2g+2(n^2-1)(g-1)$.  

 The characteristic polynomial $\det (x-\Phi)=x^n+a_1x^{n-1}+\dots +a_n=0$ defines the  {\it spectral curve} $S$ in the total space of the canonical bundle $p: K\rightarrow \Sigma$. It is the divisor of a section of $p^*K^n$ and since a   cotangent bundle  has trivial canonical bundle, it follows that  $K_S\cong p^*K^n$ and its genus is $g(S)=n^2(g-1)+1$. On $S$, by definition, 
$x$ is a single-valued eigenvalue of $\Phi$ and is the tautological section of $p^*K$ on $K$. The coefficient $a_i$ in the characteristic polynomial is a holomorphic section of $K^i$ on $\Sigma$ and $a_1,\dots,a_n$ defines a map $\pi$  from ${\mathcal M}$ to the vector space
$$B=\bigoplus_{i=1}^n H^0(\Sigma, K^i)$$
which has dimension $g+3(g-1)+\dots = g+(n^2-1)(g-1)=\dim {\mathcal M}/2$. A point of $B$ defines the equation of the spectral curve, which is generically smooth.

On the spectral curve $S$ we have an exact sequence (see \cite{BNR})
\begin{equation}
0\rightarrow U\otimes p^*K^{1-n}\rightarrow p^*E\stackrel{x-\Phi}\rightarrow p^*(E\otimes K)\rightarrow  U\otimes p^*K\rightarrow 0
\label{exact}
\end{equation}
The line bundle $U$ allows us to recover the vector bundle: $E$ is the direct image sheaf $p_*U$. From Grothendieck-Riemann-Roch, $\deg U =\deg E+(n-n^2)(1-g)$.  This describes the abelian variety for $GL(n)$ -- the direct image of any line bundle of this degree defines a stable Higgs bundle, and so the fibre of $\pi: {\mathcal M}\rightarrow B$ is isomorphic to the Jacobian of the spectral curve $S$. Functions on the base $B$ Poisson-commute and this is the description of the integrable system.

The dual $A^{\vee}$ of an abelian variety $A$ is the moduli space of degree zero holomorphic line bundles over $A$.  If $L$ is an ample line bundle on $A$, and $T_x:A\rightarrow A$ is translation by $x$, then $x\mapsto T^*_xL\otimes L^*$ identifies $ A^{\vee}$ as a quotient of $A$ by a finite subgroup. If $L$ is the  theta-divisor of the Jacobian of a curve, this map is an isomorphism so a Jacobian is its own dual.  For $G=GL(n)$, then the abelian variety is self-dual which agrees with the Langlands duality between $GL(n)$ and itself. 
\begin{rmk}
Dualizing the sequence (\ref{exact}) and tensoring with $p^*K$ gives:
\begin{equation}
0\rightarrow U^*\rightarrow p^*E^*\rightarrow p^*(E^*\otimes K)\rightarrow  U^*\otimes p^*K^n\rightarrow 0
\label{exactd}
\end{equation}
and $E^*$ is the direct image sheaf $p_* (U^*\otimes p^*K^{n-1})$.  This arises from the dual action of $\Phi$ on $E^*$.
\end{rmk}

For a  simple Lie group $G$ of rank $k$ we have, instead of the vector bundle $E$, a holomorphic principal $G$-bundle $P$ and a Higgs field $\Phi \in H^0(\Sigma,\lie{g}\otimes K)$, where  $\lie{g}$ denotes the adjoint bundle associated to $P$. The map  $\pi$ is defined by taking   $p_1,\dots, p_k$ to be a basis for the invariant polynomials  on $\lie{g}$. If $p_i$ has degree $d_i$ then evaluating on $\Phi$ we define
 
     $$\pi(\Phi)\in \bigoplus_1^n H^0(\Sigma,K^{d_i}).$$
     This is again a vector space of half the dimension $2\dim G(g-1)$ of  ${\mathcal M}$ \cite{H2}. 
     
   \begin{ex}  For the case of $G=G_2$,   $k = 2$ and $\dim G=14$, so 
     $\dim{\mathcal M}=28(g-1)$. The invariant polynomials are  
      $p_2,p_6$
    and then 
     $$\dim (H^0(\Sigma,K^2)\oplus  H^0(\Sigma,K^6))
     =3(g-1)+11(g-1)=14(g-1).$$
     \end{ex}
     
     \subsection{Prym varieties} \label{prym}
     
     For many groups, the abelian variety is related to a Prym variety, so we recall here the basic properties of these.  If $p:Y\rightarrow X$ is a degree $n$ map of compact Riemann surfaces  then there is the norm map $\Nm$ (or $\Nm_p$ when we want to keep track of the map $p$) 
  defined on divisor classes  by    
      $\Nm(\sum a_ix_i)=\sum a_ip(x_i).$  
The Prym variety $P(Y,X)$ is defined to be the connected component of the kernel of $\Nm:J(Y)\rightarrow J(X)$.  When $p^*:J(X)\rightarrow J(Y)$ is injective, which will always be the case for us (see \cite{BNR} and \cite{Lange} for the exact criteria), the Prym variety is connected.

Using  the isomorphism $J\cong J^{\vee}$, the pull-back map $p^*:J(X)\rightarrow J(Y)$ is dual to the norm map (see \cite{Lange}) and so the dual of the Prym variety is 
$ P^{\vee}(Y,X)=J(Y)/p^*J(X)$.
Restricting to $P(Y,X)\subseteq J(Y)$,   we get
$$ P^{\vee}(Y,X)=P(Y,X)/p^*J(X)\cap P(Y,X).$$
But $x\in p^*J(X)\cap P(Y,X)$ if and only if  $\Nm p^*x=0$, which, since   $\Nm (p^{-1}(x))=p(p^{-1}(x))=nx$, is when $nx=0$. In this case then, 
$ P^{\vee}(Y,X)$ is isomorphic to the quotient of  $P(Y,X)$ by the finite subgroup of elements of order $n$ in $p^*J(X)$.

\begin{ex} If $G=SL(n)$ then $\Lambda^n E$ is trivial and $\tr \Phi=0$. The abelian variety consists of line bundles $U$ on the spectral curve such that $\Lambda^n p_*U$ is trivial. But (see \cite{BNR}), for a map $p:Y\rightarrow X$ 
\begin{equation}
\Nm(U)= \Lambda^n p_*U\otimes \delta^{-1}
\label{nm}
\end{equation}
where $\delta^{-1}=\Lambda^n p_*{\mathcal O}_Y$. In the  case of $p:S\rightarrow \Sigma$ this is $K^{n(n-1)/2}$ so that  $\Lambda^n E$ is trivial if and only if
$\Nm(U)=K^{-n(n-1)/2}$, or equivalently that $U\otimes p^*K^{(n-1)/2}$ lies in the Prym variety.

The Langlands dual of $SL(n)$ is $PGL(n)$ and the dual of the Prym variety is its quotient by the elements of order $n$ in $J(\Sigma)$. But two $SL(n)$ bundles $E,E'$ are projectively equivalent if $E'=E\otimes L$ for a line bundle on $\Sigma$ of order $n$. This demonstrates the duality result for $SL(n)$.
  \end{ex}

\section{The group $Sp(2m)$ \label{symplectic}}
For the group $Sp(2m)$ we take the vector bundle $E$ to be  of rank $2m$ with a nondegenerate skew form $\langle\,\,,\,\,\rangle$ and  $\Phi$ to  satisfy $\langle \Phi v,w\rangle +\langle  v,\Phi w\rangle =0$. 

If $v_i,v_j$ are eigenvectors of $a\in\lie{sp}(2m)$,  then
$$\lambda_i\langle v_i,v_j\rangle=\langle a v,w\rangle =-\langle  v,a w\rangle=-\lambda_j \langle v_i,v_j\rangle$$
and so $\langle v_i,v_j\rangle=0$ unless $\lambda_i=-\lambda_j$. Since the skew form is nondegenerate, when the eigenvalues are distinct they must occur in opposite pairs, and so the characteristic polynomial of $\Phi$ is of the form 
$$\det (x-\Phi)=x^{2m}+a_2x^{2m-2}+\dots+a_{2m}.$$

The spectral curve $S$ defined by the above equation thus has an involution $\sigma$ defined by $\sigma(x) = -x$, and the eigenspace $L$ for $\Phi$ with eigenvalue $x$ is transformed to $\sigma^*L$ for eigenvalue $-x$. From (\ref{exact}) and (\ref{exactd}) this means that $U^*\cong U\otimes p^*K^{1-2m}$, or  $U^2\cong p^*K^{2m-1}$.
Choosing a square root $K^{1/2}$, the bundle $L_0=U\otimes p^*K^{-m+1/2}$ satisfies $\sigma^*L_0\cong L_0^*$. 

The subvariety of $J(S)$  satisfying this  condition is the Prym variety $P(S, S/\sigma)$ of the quotient map  $\pi:S\rightarrow S/\sigma$, since in this case 
$\pi^*\Nm(x)=x+\sigma x$.  For brevity we shall write $S/\sigma=\bar S$. 
\vskip .25cm
Given $L_0$ in the Prym variety we reconstruct $E$ as $p_*U$ where $U=p^*K^{m-1/2}\otimes L_0$. But since $U^2\cong p^*K^{2m-1}$  we have, from (\ref{exact}) and (\ref{exactd}), an isomorphism $E\cong E^*$ which defines the symplectic form.

\section{The group $SO(2m+1)$. \label{solie}}
\subsection{The spectral curve}
Now suppose we have a holomorphic vector  bundle $V$ of rank $2m+1$, with $\Lambda^{2m+1}V^*$ trivial and a nondegenerate symmetric bilinear form $g(v,w)$ such that $\Phi$ satisfies $g(\Phi v,w)+g(v,\Phi w)=0$. The moduli space here has two components, characterized by a class $w_2\in H^2(\Sigma,\Z_2)\cong \Z_2$, depending on whether $V$ has a lift to a spin bundle or not.

First we discuss the Lie algebra $\lie{so}(2m+1)$. Where the eigenvalues of $a\in \lie{so}(2m+1)$ are distinct, an argument like the symplectic one shows that if $a$ has distinct eigenvalues then one is zero  and the others are in opposite pairs so that the characteristic polynomial of $\Phi$ is of the form 
$$\det (x-\Phi)=x(x^{2m}+a_2x^{2m-2}+\dots+a_{2m}).$$
It is the zero eigenspace which links $\lie{so}(2m+1)$ to $\lie{sp}(2m)$ in the duality. Let $V$ be the $2m+1$-dimensional orthogonal vector space on which $SO(2m+1)$ acts and let   $V_0$ be the one-dimensional zero eigenspace of $a$, then  $a:V/V_0\rightarrow V/V_0$ is invertible and $g(av,w)$ is a non-degenerate skew form $\omega$ on the $2m$-dimensional space $V/V_0$. Since 
$$\omega(av,w)+\omega(v,aw)=g(a^2v,w)+g(av,aw)= -g(av,aw)+g(av,aw)=0$$
it follows that $a$ acts as a transformation $a'\in \lie{sp}(2m)$.

There is a canonically defined vector in $V_0$, algebraically determined by $a$. Let  $\alpha(v,w)=g(av,w)$ define the skew form $\alpha \in \Lambda^2 V^*$. Then $\alpha^m$ lies in $\Lambda^{2m} V^*$.  Let $\nu$ be the $SO(2m+1)$-invariant volume form in $\Lambda^{2m+1}V^*$. Then  $\alpha^m=m!i_{v_0}\nu$ for a unique vector $v_0$ and since   $a$ acts trivially on $\alpha$, $v_0$ is a zero eigenvector for $a$. 

Writing down $\alpha$ in an orthonormal basis $e_0,e_1,\dots e_{2m}$ gives 
$$\alpha=i\lambda_1e_1\wedge e_2+i\lambda_2 e_3\wedge e_4+\dots +i\lambda_m e_{2m-1}\wedge e_{2m}$$
where $\pm \lambda_i$ are the non-zero eigenvalues of $a$. Then  if $\nu=e_0\wedge e_1\wedge\dots \wedge e_{2m}$ 
\begin{equation}
v_0=i^m\lambda_1\lambda_2\dots\lambda_m e_{0}.
\label{v0}
\end{equation}
 In particular,  
   if $a$ has distinct eigenvalues, $\lambda_i\ne 0$ and  $v_0$ is non-null.

 Now $(V/V_0)^*\subset V^*$ is naturally the annihilator of $V_0$, and using the inner product on $V$ to identify $V^*$ with $V$ this is the orthogonal complement of $V_0$. On the other hand the symplectic form $\omega$ on $V/V_0$ identifies it with its dual. It is straightforward to see that, with this identification, the inner product restricted to $V_0^{\perp}$ can be written
$$g(u,u)=\omega(a'u,u).$$   

\vskip .15cm
We now put this into global effect for an $SO(2m+1)$ Higgs bundle $V$. In this case   $\Phi \in H^0(\Sigma, {\lie g}\otimes K)$  replaces $a$ and defines  a section $\phi$ of $\Lambda^2 V^*\otimes K$, so $\phi^m$ defines a zero eigenvector $v_0\in V\otimes K^m$.  Put another way, $v_0$ is an isomorphism from $K^{-m}$ to the zero eigenspace bundle $V_0\subset V$. 

 The global version of  $\omega$ is  $g(\Phi v,w)$, which is a skew form on $V/V_0$ with values in $K$: thus, choosing a square root of $K$, $E=V/V_0\otimes K^{-1/2}$ has a skew form which is  generically non-degenerate. But $\Lambda^{2m+1}V$ is trivial and so $\Lambda^{2m}(V/V_0)\cong V_0^*\cong K^m$, which means $\Lambda^{2m}E$ is trivial and the skew form must be non-degenerate everywhere.  We write $V_1=V/V_0$ and then
$$V_1\cong E\otimes K^{1/2}$$
where $E$ is a symplectic bundle.

As above, $\Phi$ induces a transformation $\Phi'$ on $E$ and has characteristic polynomial $x^{2m}+a_2x^{2m-2}+\dots+a_{2m}$. We are therefore precisely in the symplectic case and we can describe this structure equivalently by a line bundle $L_0$ in the Prym variety $P(S,\bar S)$ of the spectral curve $S$. The only difference  is that in Section \ref{symplectic} we chose a square root of $K$ to define $E$, but $L_0$ in the Prym variety was also defined by choosing a square root; in our case, by choosing the same square root  each time, we have a canonical $L_0$ in the Prym variety. 

\subsection{Reconstructing the bundle.}\label{reconstruct}

There  remains the task of reconstructing the bundle $V$ with $SO(2m+1)$ structure from the symplectic bundle $E$, so we examine this more closely.  Since $V_1=V/V_0$ we have an extension of vector bundles
$$0\rightarrow V_0\rightarrow V\rightarrow V_1\rightarrow 0$$
or, since $V\cong V^*$, dualizing 
$$0\rightarrow V_1^*\rightarrow V\rightarrow V_0^*\rightarrow 0.$$
 In this second picture, $V_1^*$ is the orthogonal complement of the rank one subbundle $V_0$. Thus, where  $V_0$ is non-null it splits  the sequence, which means that the extension class  is supported on the divisor   $D$ where $V_0$ is null. From (\ref{v0}) this  has the equation $a_{2m}=\lambda_1^2\dots\lambda_m^2=0$.

Over $D$, $V_0$ is null and therefore contained in its orthogonal complement, i.e. we have an inclusion $V_0\subset V_1^*$, which we regard as a section $i \in H^0(D,\Hom(V_0,V_1^*))$. We also view $a_{2m}\in H^0(\Sigma,K^{2m})$ as a homomorphism from $V_0 \cong K^{-m}$ to $V_0^*\cong K^m$. Consider now the exact sequence of sheaves   
$$0\rightarrow {\mathcal O}_{\Sigma}(\Hom(V_0^*,V_1^*))\stackrel{a_{2m}}\rightarrow {\mathcal O}_{\Sigma}(\Hom(V_0,V_1^*))\rightarrow {\mathcal O}_D(\Hom(V_0,V_1^*))\rightarrow 0.$$
In the long exact cohomology sequence 
\begin{equation}
\rightarrow H^0(D, \Hom(V_0,V_1^*))\rightarrow H^1({\Sigma},\Hom(V_0^*,V_1^*))\stackrel{a_{2m}}\rightarrow H^1({\Sigma}, \Hom(V_0,V_1^*))\rightarrow 
\label{seq}
\end{equation}
the section $i \in H^0(D,\Hom(V_0,V_1^*))$ defines a class $\delta(i)\in H^1({\Sigma},\Hom(V_0^*,V_1^*))$ which we claim defines the extension $V$.  To see this,  cover  ${\Sigma}$ by a union $U$ of small discs centred on the points of $D$, together with the single open set ${\Sigma}\setminus D$.  As remarked above, the sequence is split on ${\Sigma}\setminus D$. Choose a splitting over $U$, then over $U$ the inclusion $V_0\subset V$ can be written as $s\mapsto (u(s),a_{2m}s)\in V_1^*\oplus V_0^*$. Here, $u$ is a holomorphic extension from $D$ to $U$ of the inclusion $i$. Since $V_0$ defines the splitting outside $D$, a \v Cech cocycle for the extension is defined by 
$$u(s)/a_{2m}\in H^0(U\cap  {\Sigma}\setminus D,  \Hom(V_0,V_1^*))$$
and this is a representative for $\delta(i)$.
\vskip .25cm
 We have just seen that to construct  the orthogonal bundle we   use a homomorphism $K^{-m}\rightarrow V_1^*\cong E\otimes K^{-1/2}$ on $D$. It is this that we need to focus on, and identify from the symplectic viewpoint.

If $a_{2m}$ has simple zeros (which will be so if  the spectral curve is smooth), then at a point on $D$, the $SO(2m+1)$ Higgs field $\Phi$ decomposes $V$ into a direct sum of orthogonal invariant subspaces $V=U_0\oplus U_2\oplus \dots \oplus U_{m}$ where $\dim U_0=3$ and $\Phi$ restricted to $U_0$ is nilpotent, and where $U_i$, for $i>0$,  is the sum of  the $\pm \lambda_i$ eigenspaces. There is an orthonormal basis  ${e}_0,{e}_1,{e}_2$  for $U_0$ such that, using the usual three-dimensional  vector cross product,  
$$\Phi({ x})=\mu ({ e}_1+i{ e}_2)\times { x}.$$
Then the two-form defined by $\Phi$ is 
$$\phi=i\mu e_0\wedge (e_1+ie_2)+i\lambda_2 e_3\wedge e_4+\dots +i\lambda_m e_{2m-1}\wedge e_{2m}$$
and we obtain
\begin{equation}
v_0=-i^{m-1}\mu\lambda_2\lambda_3\dots\lambda_m (e_1+ie_2)
\label{v00}
\end{equation}
Now consider the two-dimensional space  $V_0^{\perp}\cap U_0$. This is spanned by $e_1+ie_2$ and $e_0$. The vector $u=i^m\lambda_2\lambda_3\dots\lambda_m e_0$ has the property $\Phi(u)=v_0$ and 
$(u,u)=(-1)^m\lambda_2^2\dots\lambda_m^2$. But note that where $a_{2m}=0$ (say $\lambda_1=0$), $\lambda_2^2\dots\lambda_m^2=(-1)^{m-1}a_{2m-2}$, using the   coefficient $a_{2m-2}$ in the characteristic polynomial. The vector $u$ is  defined only modulo $e_1+ie_2$ by these properties, so there is a distinguished non-zero vector $u$ in the one-dimensional space $(V_0^{\perp}\cap U_0/V_0)\otimes K^{m-1}$ such that $\Phi(u)=v_0\in V_0\otimes K^m$ and $(u,u)=-a_{2m-2}$.

This data is visible from the symplectic viewpoint, but not quite uniquely determined.  By definition $V_0^{\perp}=E\otimes K^{-1/2}$ and we have at a point of $D$ a symplectic-orthogonal $\Phi'$-invariant decomposition $E=E_0\oplus E_2\oplus \dots E_m$, where $\Phi'$ is nilpotent on $E_0$. Here we look for an $e\in E\otimes K^{m-3/2}$ such that $\omega(\Phi'(e),e)= - a_{2m-2}$. There are two possible choices  and since $D$ (a divisor of $K^{2m}$) has degree $4m(g-1)$ there is a total of $2^{4m(g-1)}$ such choices. Each one under $\Phi$ defines a $v_0\in V_0\otimes K^m$ to construct an orthogonal bundle. Here is the essential point, missed out in  \cite{H2}, and which in the next section we shall see describes the duality. 

\begin{thm} \label{canon} Let $(E,\Phi')$ be a generic symplectic Higgs bundle  of rank $2m$. Then an associated $SO(2m+1)$ Higgs bundle is determined by a vector  $e\in E_a\otimes K_a^{m-3/2}$, for each point $a\in \Sigma$ where $a_{2m}(a)=0$, such that $\omega_a(\Phi'(e),e)= - a_{2m-2}(a)$ 
\end{thm}
\begin{prf}
 We have constructed an extension $V$ 
 \begin{equation}
 0\rightarrow V_1^*\rightarrow V\rightarrow K^m\rightarrow 0
 \label{ext0}
 \end{equation} 
 with this information. We need to define  a metric and a Higgs field  on  $V$.

First the metric. 
In the exact cohomology sequence (\ref{seq})  $\delta(i)\in H^1(\Sigma, \Hom(K^m,V_1^*))$  vanishes when multiplied by $a_{2m}$. This means that we can lift the homomorphism $a_{2m}:K^{-m}\rightarrow K^{m}$  to a  homomorphism $\alpha: K^{-m}\rightarrow V$. There are many of these -- any two liftings will differ by a homomorphism from $K^{-m}$ to $V_1^*$ --  but we shall choose a distinguished one later, after constructing the metric.

If we use the concrete cocycle description of the extension, with the local splitting $V\cong V_1^*\oplus K^m$  near $D$, then $\alpha$ has the form
\begin{equation} 
\alpha(s)=\left(a_{2m}\frac{u(s)}{a_{2m}}, a_{2m}s\right)=(u(s),a_{2m}s)
\label{inc}
\end{equation}
and so maps $K^{-m}$ isomorphically to a rank one subbundle $V_0$. 

Outside of $D$, we now define, using the symplectic structure on $E$,  an inner product on $V$ so that $V_1^*$ and $V_0$ are orthogonal:  for $v\in V_1^*, s\in K^{-m}$, put  
\begin{equation}
g(v+s,v+s)=\omega(\Phi' v,v) + (-1)^m a_{2m}s^2.
\label{metric}
\end{equation}
We must show that this extends over $D$.

We use a local coordinate $z$ in a neighbourhood of a point $z=0$ of $D$ and  trivialize $K$ with $dz$. We shall prove local regularity by using the local splitting $V_1^*\oplus K^m$ near $D$. We can take $v=0,s=1$ as the canonical zero eigenvector $v_0=\phi_0e$ defined above, since in the metric (\ref{metric}) $g(v_0,v_0)=(-1)^ma_{2m}$. This means that we can be more explicit about the particular form (\ref{inc}) of the inclusion of $K^{-m}$: there is a local nonvanishing section of the line bundle $V_0$ of the form
$$(v_0+zv_{01}+\dots, a_{2m}(z))=(v_0,0)+z(v_{01},c)+\dots$$
where $a_{2m}=cz+\dots$ and $c\ne 0$.

A vector $(w,t)$ in the splitting $V_1^*\oplus K^m$ can thus be written in the orthogonal splitting $V_1^*\oplus V_0$ as $(v,s)$ where
$$v= w - \frac{t}{cz}(v_0+zv_{01}+\dots),\qquad s=\frac{t}{cz}+\dots$$

The inner product  (\ref{metric}) evaluated on $(v,t)\in V_1^*\oplus K^m$ is now
$$\omega\left(\Phi' (w - \frac{t}{cz}(v_0+zv_{01}+\dots),w - \frac{t}{cz}(v_0+zv_{01}+\dots)\right) +(-1)^m (cz+\dots)\left(\frac{t^2}{c^2z^2}+\dots\right)$$
and we have to show that, despite the denominators $z$, this is smooth.

Near $z=0$ we have 
$ \Phi'=\phi_0+z\phi_1+\dots$
where $\phi_0$ has a one-dimensional kernel spanned by $v_0$. Thus
$$ \Phi' (v_0+zv_{01}+\dots) =z(\phi_1 v_0+\phi_0 v_{01})+\dots$$
so all we need to show is that
$$\frac{1}{c^2z^{2}}\omega( z(\phi_1 v_0+\phi_0 v_{01}), v_0+ zv_{01}) +(-1)^m\frac{1}{ cz}$$
is smooth. But
$$\omega(\phi_0 v_{01},v_0)=-\omega(v_{01},\phi_0 v_0)=0$$
so for regularity we just need to show that 
\begin{equation}
\frac{1}{c^2 }\omega(\phi_1 v_0, v_0) +(-1)^m\frac{1}{ c}
\label{vanish}
\end{equation}
vanishes.

Now consider $\det \Phi'= a_{2m}=cz+\dots$. Choose a basis where $e_2=e, e_1=\phi_0 e$ and the others are eigenvectors. We then have 
\begin{eqnarray*}
\det  \Phi'(e_1\wedge\dots\wedge e_{2m})&=& \Phi' e_1\wedge  \dots \wedge  \Phi' e_{2m}\\
&=&(-1)^{m-1}z\lambda_2^2\dots\lambda_m^2 \phi_1 e_1\wedge \phi_0 e_2\wedge e_3\dots \wedge e_{2m}+\dots\\
&=& (-1)^{m-1}z\lambda_2^2\dots\lambda_m^2\phi_1e_1\wedge e_1\wedge e_3\wedge\dots e_{2m}\\
&=&(-1)^{m-1}z\lambda_2^2\dots\lambda_m^2\frac{\omega(\phi_1 e_1,e_1)}{\omega(e_1,e_2)}e_1\wedge e_2\wedge e_3\dots\wedge e_{2m}
\end{eqnarray*}
But $\omega(e_1,e_2)=\omega(\phi_0 e, e)=(-1)^m\lambda_2^2\dots \lambda_m^2$ by the choice of $e$ and so $\omega(\phi_1v_0,v_0)=(-1)^{m-1} c$ and the term (\ref{vanish}) does indeed vanish.

To complete the construction of the metric, observe that the bilinear form defines a homomorphism from $V$ to $V^*$ but 
 $V$ satisfies $\Lambda^{2m+1}V\cong \Lambda^{2m}V_1^*\otimes K^m$, which is trivial since $\Lambda^{2m}V_1^*\cong K^{-m}$. Thus this homomorphism has everywhere non-zero determinant and the form is non-singular everywhere. 

Now we define the Higgs field. The metric identifies $V$ with $V^*$ so we have a dual extension  to (\ref{ext0}):
$$0\rightarrow V_0\rightarrow V\rightarrow V_1\rightarrow 0.$$

We have the symplectic Higgs field $\Phi':V_1\rightarrow V_1\otimes K$ and we let 
$\Phi: V\rightarrow V\otimes K$ be the composition
$$V\rightarrow V/V_0=V_1\rightarrow V_1\otimes K \rightarrow V_1^*\otimes K\subset V\otimes K$$
where the second arrow is $\Phi'$ and the third is the inner product $\omega(\Phi'v,v)$ on $V_1^*=E\otimes K^{-1/2}$.
\vskip .25cm
We have seen here how any lifting $\alpha:K^{-m}\rightarrow V$ leads to a metric and a Higgs field. We shall now see that there is a unique lift such that  the composition $\Phi\alpha=0$. 

Take any lift and set $\Psi=\Phi\alpha$. Since $\Phi(V)\subseteq V_1^*\otimes K$, $\Psi\in H^0(\Sigma, \Hom(K^{-m}, V_1^*\otimes K))$ which is  $H^0(\Sigma,  V_1^*\otimes K^{m+1})$. Another lifting differs from $\alpha$ by $\beta\in H^0(\Sigma,  V_1^*\otimes K^{m})$ so there exists $\beta$ with $\Phi(\alpha-\beta)=0$ if $\Psi$ is in the image of $\Phi':H^0(\Sigma,V_1^*\otimes K^m)\rightarrow H^0(\Sigma, V_1^*\otimes K^{m+1})$. Consider the sequence of sheaves 
$$0\rightarrow {\mathcal O}_{\Sigma}(V_1^*\otimes K^m)\stackrel{\Phi'}\rightarrow {\mathcal O}_{\Sigma}(V_1^*\otimes K^{m+1})\rightarrow {\mathcal O}_D(S)\rightarrow 0$$
where $S$ is the skyscraper sheaf of cokernels of $\Phi'$ at $D$ (recall that $a_{2m}=0$ is precisely where $\Phi'$ has a zero eigenvalue.) By construction, on $D$, $\alpha$ takes values in the kernel of $\Phi'$, so from the exact cohomology sequence there is a unique $\beta$ for which $\Phi (\alpha-\beta)=0$.
\end{prf}

\subsection{Duality}\label{duality}

We shall now show  that the data for an $SO(2m+1)$ Higgs bundle above is given by a point in the \emph{dual} of the Prym variety $P(S,\bar S)$ for an $Sp(2m)$ bundle, thus giving a realization of Langlands duality within this context. It was Michael Thaddeus \cite{Th} who pointed out a mistake on page 108 of the author's paper \cite{H2}, the resolution of which  yields duality of the abelian varieties concerned and not their equality as stated on page 109 of that paper. 

From Theorem \ref{canon} the extra data for constructing an $SO(2m+1)$ bundle from a symplectic bundle is a choice between two  vectors $\pm e$ at each point of the divisor $D$. 

The symplectic bundle $E$ was defined as $V_1^*\otimes K^{1/2}$ and  this eigenspace bundle, pulled back to $S$, is $U\otimes p^*K^{1-2m}$ where  $U=p^*K^{m-1/2}\otimes L_0$ and $L_0$ lies in the Prym variety $P=P(S,\bar S)$. 

We can identify via the projection $p:S\rightarrow \Sigma$ the finite set of points $D$ on $\Sigma$ defined by $a_{2m}=0$ with the zero section $x=0$ on the spectral curve $S$. To avoid confusion we shall call this the divisor $D_S$ (of $p^*K$). Then  we see that on $D_S$ there is a natural isomorphism
of $V_0$ with $K^{-m}\otimes L_0$. This means 
 that our choice of isomorphism $V_0\cong  K^{-m}$ is the same as a choice of trivialization of $L_0$ on $D_S$.
The trivialization is not arbitrary -- it satisfies the quadratic condition $\omega(\Phi'(e),e)= - a_{2m-2}$ given in  Theorem \ref{canon}. 

Now $D_S$, defined by $a_{2m}=0$  is  the fixed point set  of  $\sigma(x)=-x$ on the spectral curve $x^{2m}+a_2x^{2m-2}+\dots+a_{2m}=0$. For a  line bundle on $S$ in the Prym variety there is by definition an isomorphism from $\sigma^*L$ to $L^*$ and so at the fixed points we have  an  isomorphism $L\cong L^*$, or equivalently a non-zero section $u_L$ of $L^2$ on $D_S$. The quadratic condition is that  we have to trivialize $L_0$ by choosing a section  $v$ of $L_0$ on $D_S$ such that $v^2=u_{L_0}$.
\vskip .25cm
This data, a point $L_0\in P(S,\bar S)$ and a trivialization of $L_0$ over $D_S$, defines a finite covering of $P(S,\bar S)$ of degree  $2^{4m(g-1)}$. It is also a group under tensor product and the covering is a homomorphism. Now a trivialization on $D_S$ multiplied by $-1$ gives a scalar multiple of the extension class $\delta(i)$ of Section \ref{reconstruct} and hence the same vector bundle $V$, so the data for constructing $V$  actually lies in a covering $P'$ of degree $2^{4m(g-1)-1}$

There is one straightforward way to find elements in $P'$:  if $L_0=N^2$ for some line bundle $N\in P(S,\bar S)$, then we can take $v=u_N$, so that $v^2=u_{N^2}=u_{L_0}$. So consider the squaring map $s:P(S,\bar S)\rightarrow P(S,\bar S)$ defined by $s(L)=L^2$. This is surjective since the Prym variety is connected. Its kernel consists of equivalence classes of line bundles for which $\sigma^*L\cong L^*$ and $L^2$ is trivial. The latter condition is $L^*\cong L$ and together with the first we obtain an isomorphism 
$$\sigma^*L\cong L.$$
This defines a lifting $\tilde \sigma$ of the action of $\sigma$ on the curve $S$ to the line bundle $L$. The trivialization of $L^2\cong {\mathcal O}$ at the fixed point set $D_S$ of $\sigma$ is then just the action $\pm 1$ of $\tilde \sigma$.  But if the action is trivial at all points of $D_S$, the line bundle $L$ is pulled back from $\bar S$. It follows that the quotient
$$P(S,\bar S)/\pi^*H^1(\bar S,\Z_2)$$
maps injectively to $P'$. 
\vskip .25cm
From the Riemann-Hurwitz formula for the covering $S\rightarrow \bar S$ the genus of $\bar S$ is given by 
$$2g(\bar S)=g(S)+1-2m(g-1)=(4m^2-2m)(g-1)+2$$
and the dimension of the Prym variety is
$$g(S)-g(\bar S)=4m^2(g-1)+1-(2m^2-m)(g-1)-1=2m(m+1)(g-1)$$
(which is of course $\dim Sp(2m)(g-1)$). Thus $P(S,\bar S)/\pi^*H^1(\bar S,\Z_2)$ projects under the squaring map to $P(S,\bar S)$ as a covering of degree $$2^{2(g(S)-2g(\bar S))}=2^{4m(g-1)-2}.$$
This is half of the degree of the covering $P'$. The reason is that $P'$ has two components -- the data that gives a spin bundle and its complement. Since $P(S,\bar S)$ is connected, its image has constant $w_2$. However, either component is acted on freely and transitively by  $P(S,\bar S)/\pi^*H^1(\bar S,\Z_2)$, and this is, as  we saw in \ref{prym},  the {\it dual} of $P(S,\bar S)$. 

Thus finally we see how the duality of abelian varieties corresponds to Langlands duality for $Sp(2m)$ and $SO(2m+1)$.

\begin{rmks}  
\begin{enumerate}
\item It is in  fact the identity component of $P'$ which corresponds to spin bundles.   The natural origin of the Prym variety is a point in the Teichm\"uller component of \cite{HLie},  since $SO(2m+1)$ is the adjoint group. The vector bundle is 
$$V=K^{-m}\oplus K^{-m+1}\oplus \dots\oplus 1\oplus \dots \oplus K^{m-1}\oplus K^m$$
(with the obvious pairings defining the metric) and the Higgs field is a canonical normal form for the given characteristic polynomial. The point to notice here is that 
$$V=1\oplus (K\oplus K^{-1})\oplus (K^2\oplus K^{-2})\oplus \dots$$  
is an orthogonal sum of  $SO(2)$ bundles $K^n\oplus K^{-n}$, each of which is spin, indeed $K^{\pm n/2}$ are the two spin bundles. So $w_2=0$.
\item
The two components are covering spaces whose group of order $2^{4m(g-1)-2}$ consists of the elements in $H^0(D,\Z_2)$ with an even number of minus signs, modulo the constant functions $\Z_2$. This follows from the interpretation as the action, at the fixed point set, of $\tilde\sigma$ on a line bundle $L$ of degree zero.  If $n_+,n_-$ are the numbers of points of $D$ with action $+1,-1$ respectively then the Lefschetz fixed point formula gives
$$\frac{1}{2}(n_+-n_-)=\tr \tilde \sigma\vert_{H^0(L)}-\tr \tilde \sigma\vert_{H^1(L)}=N_+-N_-$$ 
where
$$n_++n_-=4m(g-1)\qquad N_++N_-=1-g(S)=-4m^2(g-1).$$
Hence $n_-=(4m^2+2m)(g-1)+2N_-$ is {\it even}.
\end{enumerate}
\end{rmks}

\section{The group $G_2$}
\subsection{The geometry of $G_2$ \label {gee2}}
In \cite{Katz}, Katzarkov and Pantev gave one description of the abelian variety  which defines a $G_2$ Higgs bundle. We shall achieve the same end, but use less Lie theory. Our point of view will be that in many respects $G_2$ is not an exceptional Lie group, and dealing with it head-on as in the case of the classical groups, we shall be able to see more closely what is happening.

Our starting point is that the complex group $G_2$ is the connected component of the subgroup of $GL(7)$ which preserves a generic $3$-form $\rho$ on $\C^7$ -- in other words $\rho$ lies in an open orbit in the space of all three-forms (see for example \cite{H4},\cite{KS}). The form defines a metric on a $7$-dimensional vector space $V$ as follows.

If $v\in V$ then $i_v\rho\wedge i_v\rho\wedge\rho\in \Lambda^7V^*$. This is a quadratic form $c(v,v)$ with values in $\Lambda^7V^*$ and so defines a map $V\rightarrow V^*\otimes \Lambda^7V^*$ whose determinant lies in $(\Lambda^7V^*)^9$. This equivariant polynomial in $\rho$ has degree $21$ but is in fact the third power of a polynomial $\kappa(\rho)$ of degree $7$. The metric is defined by 
$g=c/\kappa^{1/3}$. The stabilizer of $\rho$ is  the group $G_2\times \Z_3$ with $\Z_3$ acting non-trivially on the cube root of $\kappa$. The connected component $G_2$ preserves the metric, and by construction a volume form and thus lies in $SO(7)$. 
 The open orbit in $\Lambda^3V^*$ is defined by $\kappa(\rho)\ne 0$. 
\vskip .15cm

For three-forms in six dimensions there is a  similar story -- for a $6$-dimensional complex vector space $W$, there is an open orbit in $\Lambda^3 W^*$  under the action of $GL(6)$ whose stabilizer is $SL(3)\times SL(3)\times \Z_2$ (see \cite{H3},\cite{H4},\cite{KS},\cite{chan},\cite{Bry}).  If $x_1,x_2,x_3,y_1,y_2,y_3$ is a basis with  dual basis $\xi_1,\xi_2,\xi_3,\eta_1,\eta_2,\eta_3$, a  normal form is 
\begin{equation}
\Omega=\xi_1\wedge \xi_2\wedge \xi_3+\eta_1\wedge \eta_2\wedge \eta_3.
\label{norm1}
\end{equation}

As in \cite{H3} we define for a general three-form $\Omega$ the  linear transformation ${\mathrm K}_{\Omega}$  by
$${\KK}_{\Omega}(w)=i_w\Omega\wedge \Omega\in \Lambda^5W^*\cong W\otimes \Lambda^6W^*$$ 
and then $\KK_{\Omega}^2=\lambda(\Omega)1$, where $\lambda(\Omega)\in (\Lambda^6 W^*)^2$ is an equivariant quartic polynomial.  The open orbit in $\Lambda^3W^*$ is defined by $\lambda(\Omega)\ne 0$. 

On the hypersurface $\lambda(\Omega)=0$  there is (in the induced topology) also an open orbit with normal form
\begin{equation}
\Omega=\xi_1\wedge \eta_2\wedge \eta_3+\xi_2\wedge \eta_3\wedge \eta_1+\xi_3\wedge \eta_1\wedge \eta_2.
\label{norm2}
\end{equation}
When $\lambda(\Omega)\ne 0$, $\KK_{\Omega}$ has two three-dimensional eigenspaces on which $\Omega$  restricts to a non-vanishing form (in (\ref{norm1})  they are spanned by $x_1,x_2,x_3$ and $y_1,y_2,y_3$ respectively).  When $\lambda(\Omega)=0$, 
 $\KK^2_{\Omega}=0$ and on its open orbit in the hypersurface has three-dimensional kernel spanned by $x_1,x_2,x_3$ (from (\ref{norm2})).
\vskip .15cm

The two structures are linked. As is well-known, the compact group $G_2$ acts transitively on the sphere $S^6$ with stabilizer $SU(3)$, so the orthogonal complement of a unit vector $e_7\in \R^7$ has structure group $SU(3)$. In fact (see \cite{Sal}) the $G_2$ three-form  $\rho$ can be written as 
\begin{equation}
\rho=\Omega+\varphi\wedge e_7
\label{decomp}
\end{equation}
where $\Omega$ is the real part of the holomorphic three-form on $\C^3$ fixed by $SU(3)$ and $\varphi$ is the hermitian $2$-form.  
\vskip .25cm
We are concerned with the complexification of this picture. If we replace $e_7$ by  a non-null vector $v$ in $\C^7$ then the restriction of $\rho$ to the  orthogonal complement of $v$ is a  three-form $\Omega$ and $i_v\rho$ restricts to a  two-form $\varphi$. Note that this is not our skew form $\omega$: it is the ``hermitian"  form 
$$\varphi(u,v)=g(Iu,v)=\omega(\Phi'Iu,v)$$
where 
$I=\KK_{\Omega}/\sqrt{-\lambda(\Omega)}.$
In our case, the form $\varphi$ becomes degenerate where $v$ is null, but our $\omega$ is always symplectic. 

Under the action of the symplectic group $Sp(6)$ defined by  $\omega$,   $\Omega$ lies in an open orbit of $\C^*\times Sp(6)$ on the  $14$-dimensional space of primitive $3$-forms (i.e. $\Omega \wedge \omega=0$). Its stabilizer is $SL(3)\times \Z_2$ (see \cite{KS},\cite{Ban}). In either normal form above, we can take $\omega=\xi_1\wedge\eta_1+\xi_2\wedge\eta_2+\xi_3\wedge\eta_3$ and the eigenspaces of $\KK_{\Omega}$ are then Lagrangian.

\subsection{The Lie algebra of $G_2$ \label {geelie}}
Suppose $a$ is in the Lie algebra  $\lie{g}_2\subset \lie{so}(7)$, with distinct eigenvalues. Then as in Section \ref{solie}  it has a non-null zero eigenvector and acts on its orthogonal complement $W$ preserving the symplectic form $\varphi$. It also preserves  the three-dimensional  eigenspaces $W^+$ and $W^-$ of $\KK_{\Omega}$. Its eigenvalues on $W^+$ are $\lambda_1,\lambda_2,\lambda_3$ which satisfy $\lambda_1+\lambda_2+\lambda_3=0$ since $\Omega$ restricts to an invariant volume form there, and on $W^-$ (which is dual to $W^+$), it has eigenvalues $-\lambda_1,-\lambda_2,-\lambda_3$.

Consider the two basic invariant polynomials $$f=\lambda_1^2+\lambda_2^2+\lambda_3^2,\quad q=(\lambda_1\lambda_2\lambda_3)^2.$$ Then the characteristic polynomial of $a$ is 
\begin{equation}
x(x^6-fx^4+\frac{f^2}{4}x^2-q)
\label{charg2}
\end{equation}

\subsection{The spectral curve}

Following the previous discussion, we consider a $G_2$ Higgs bundle as a rank $7$ vector bundle $V$ with $\Lambda^7 V$ trivial and with a section $\rho$ of $\Lambda^3 V^*$ which lies in the open orbit $\kappa\ne 0$ of $GL(7)$ at each point.  Because this defines an $SO(7)$ structure, we can follow the procedures of  Section \ref{solie}, and consider the spectral curve $S$ which is a divisor in the total space of $p:K\rightarrow \Sigma$. From (\ref{charg2}) its equation is 
\begin{equation}
x^6-fx^4+\frac{f^2}{4}x^2-q=0.
\label{x}
\end{equation}
 From (\ref{standard}) it has genus $g(S)=36(g-1)+1$ and is a $6$-fold cover of $\Sigma$. 

We  define as in Section \ref{solie}  the kernel $V_0$ of $\Phi$ and the symplectic bundle $E=V_1\otimes K^{-1/2}$ with induced Higgs field $\Phi'$. We now have the extra data induced by the three-form $\Omega$, which lies in 
$ H^0(\Sigma,K^{3/2}\otimes \Lambda^3 E^*).$ We then obtain
\begin{equation}
\KK_{\Omega}:E\rightarrow E\otimes K^3.
\label{KK}
\end{equation}

\subsection{The intermediate curve}
Equation (\ref{KK}), defines a ``Higgs field" on $\Sigma$ but with $K$ replaced by $K^3$. Since $\KK^2_{\Omega}=\lambda(\Omega)1$, we have $\lambda(\Omega)\in H^0(\Sigma,K^6)$ which, as remarked above,  vanishes on $D$ and so  $\lambda$ is a multiple of the coefficient $a_6=-q$. 
\vskip .15cm

The equation $z^2=q$  defines in the total space of $K^3 \rightarrow \Sigma$ a spectral curve $C$ for $\KK_{\Omega}$, which is a double covering  of $\Sigma$  on which $\sqrt{\lambda(\Omega)}$ is well-defined. Let $p_C:C\rightarrow \Sigma$ denote the projection. On $C$ we have well-defined   rank $3$  vector bundles $W^+,W^- \subset p_C^*E$  which are eigenspaces of $\KK_{\Omega}$. 

The canonical bundle of the total space of $K^3$ is the pull-back of $K^{-2}$  so that since $C$ is the divisor of a section of $K^6$ pulled back,
$$K_{C}\cong p_C^*K^4.$$
In particular, it follows by adjunction that the genus of $C$ is $g(C)=8(g-1)+1$. The set of points $x=0$ on $C$ maps isomorphically to the divisor $D$ on $\Sigma$ but we shall call it $D_C$ on $C$. It is a divisor of $p_C^*K^3$.

We need to consider the restriction of $\Omega\in  H^0(\Sigma,K^{3/2}\otimes \Lambda^3 E^*)$ to $W^+$. Let $w$ be a local coordinate in a neighbourhood of a point of $D_C$. For a generic Higgs bundle,  $\Omega$ at $w=0$  lies in the open orbit of the hypersurface. Pull back  to $C$ and it  is of the local form
$\Omega_0+w^2\Omega_1+\dots$ where $\Omega_0$ has the normal form (\ref{norm2}), and $\ker \KK_{\Omega}$ is spanned by $x_1,x_2,x_3$. 

Let $\tilde x_i=x_i+wv_i+\dots$ be a local basis of sections for  $W^+$. Restricting $\Omega$ gives
$$\Omega(\tilde x_1,\tilde x_2,\tilde x_3)=w[\Omega_0(x_1,x_2,v_3)+\Omega_0(x_2,x_3,v_1)+\Omega_0(x_3,x_1,v_2)]+O(w^2)$$
But  the explicit normal form  (\ref{norm2}) is
$\Omega=\xi_1\wedge \eta_2\wedge \eta_3+\xi_2\wedge \eta_3\wedge \eta_1+\xi_3\wedge \eta_1\wedge \eta_2$
and so the coefficient of $w$ vanishes. Hence we have a section of $\Lambda^3 (W^+)^*\otimes p_C^*K^{3/2}$ which vanishes on $D_C$ with multiplicity $2$. Since $D_C$ is a divisor of  $p_C^*K^3$, it follows that 
\begin{equation}
\Lambda^3 W^+\cong p_C^*K^{-9/2}
\label{detw}
\end{equation}

Now $\Phi'$ preserves $W^+$ and so now we have a  ``Higgs field"  $$\Phi'':W^+\rightarrow W^+\otimes p_C^*K$$ on $C$. Its eigenvalues are eigenvalues of $\Phi'$ and indeed,
substituting $z^2=q$ in the equation of the spectral curve $S$ we have 
$$0=x^6-fx^4+\frac{f^2}{4}x^2-z^2=x^2(x^2-{f}/{2})^2-z^2$$
and 
\begin{equation}
z=x(x^2-{f}/{2}).
\label{piequation}
\end{equation}
This is an explicit degree $3$ map $p_S:S\rightarrow C$, and writing it as 
$$x^3-fx/2-z=0$$
this represents $S$ as the spectral curve of  $\Phi''$ on $C$.  

The projection $p$ from the spectral curve $S$ to $\Sigma$ therefore admits a factorization
$$S\stackrel{p_S}\rightarrow C \stackrel{p_C}\rightarrow \Sigma$$
and consequently the bundle 
$$E=(p_Cp_S)_*U=p_{C\,*}p_{S\,*}U$$
where
$p_{S\,*}U$ is a rank three vector bundle on $C$.

Since $E=p_{C\,*}p_{S\,*}U$, over $C$ there is a natural surjective homomorphism $p_C^*E\rightarrow p_{S\,*}U$ and the kernel of this is the eigenspace $W^+$ of $\KK_{\Omega}$. Thus, from  (\ref{detw})
$$\Lambda^3 p_{S\,*}U\cong \Lambda^3 (W^+)^*\cong p_C^*K^{9/2}.$$
But $p_S:S\rightarrow C$ is the spectral curve of $\Phi''\in H^0(C,\End W^+\otimes p_C^*K)$ and so
$$\Lambda^3 p_{S\,*}U=\Nm_{p_S}(U)\otimes p_C^*K^{-3}.$$
It follows that
\begin{equation}
\Nm_{p_S}(U)\cong p_C^*K^{15/2}.
\label{nmr}
\end{equation}
From the $SO(7)$ point of view, we defined $L_0=U\otimes p^*K^{-5/2}$ where $L_0\in P(S,\bar S)$ so we see from (\ref{nmr}) that 
$$\Nm_{p_S}(L_0)=\Nm_{p_S}(U\otimes p_S^*p_C^*K^{-5/2})=\Nm_{p_S}(U)\otimes p_C^*K^{-15/2}\cong {\mathcal O}.$$
This means that $L_0$ lies in the Prym variety $P(S,C)$ of $p_S:S\rightarrow C$ as well as the Prym $P(S,\bar S)$: equivalently it is the subgroup $P(S,C)^-$ -- the line bundles in $P(S,C)$ for which $\sigma^*L\cong L^*$. 

\subsection{Reconstructing the bundle}

We shall show eventually that $P(S,C)^-$ is connected and is the abelian variety for the $G_2$ Higgs bundle, but   we need now to understand the covering in order  to reconstruct the $SO(7)$ bundle. As in Section \ref{solie}, this involves the behaviour on the divisor $D$ where $a_{6}=0$, and especially the geometry of the form $\Omega$ at these points. 

\vskip .25cm
When $\Omega$ is in the singular normal form (\ref{norm2}), $\KK_{\Omega}$ has a three-dimensional kernel $U$ and $u\mapsto i_u\Omega$ gives an isomorphism 
$U\cong \Lambda^2(W/U)^*.$
But $U$ is Lagrangian so $W/U\cong U^*$ and  we get an isomorphism 
$$\ast: U\cong \Lambda^2 U.$$
Using $\Lambda^2U\cong U^*\otimes \Lambda^3U$, this defines, as in the $G_2$ argument above, a quadratic form $c(u,u)$ with values in $\Lambda^3U$, but now its determinant $\kappa$ lies in $\Lambda^3U$ so $c/\kappa$ is an inner product and $\kappa^{-1}$ a volume form on $U$. Thus $U$ acquires the standard structure of three-dimensional Euclidean space where  $\ast$ is just the Hodge star operator. 

 Choose a complementary Lagrangian subspace to $U$ and use the inner product on $U$,  then we can write $W=U\oplus U$ where the symplectic form is 
$$\omega((x_1,y_1),(x_2,y_2))=(x_1,y_2)-(x_2,y_1).$$
 The stabilizer in $GL(U)$ of the  $\ast$-operator  is $SO(3)$. Let $G$ be the stabilizer in $Sp(6)$ of $\Omega$ in this normal form. Then we have a homomorphism $G\rightarrow SO(3)$ whose kernel is of the form $(x,y)\mapsto (x+My,y)$.  To preserve $\omega$, $M$ must be symmetric. To preserve $\Omega$ in (\ref{norm2}) $M$ must have trace zero.
 Thus $G$ is the semi-direct product of $SO(3)$ with the trace-zero $3\times 3$ symmetric matrices. Its Lie algebra  consists of transformations  of the  form
\begin{equation}
({ x},{ y})\mapsto({ a
}\times { x}+M { y}, { a
}\times { y})
\label{fiform}
\end{equation}
using the vector cross product in $\C^3$.

\begin{rmk} Note that $\dim G=8$ and hence $Sp(6)$ has a  $21-8=13$-dimensional orbit passing through $\Omega$. This is the open orbit in the hypersurface $\lambda(\Omega)=0$ in the 14-dimensional space of primitive three-forms.
\end{rmk}

\vskip .25cm
Now consider the inclusion  of the zero eigenspace $V_0\cong K^{-3}$ of $\Phi$. It is defined by $v_0\in H^0(\Sigma, V\otimes K^3)$, and we can then form 
$$i_{v_0}\rho\in H^0(\Sigma, \Lambda^2V^*\otimes K^3).$$
Restrict $\rho$ to  $V_1^*=E\otimes K^{-1/2}$, the orthogonal complement to  $V_0$, and we get a form $\Omega\in H^0(\Sigma, \Lambda^3 E^*\otimes K^{3/2})$. Restrict $i_{v_0}\rho$ and 
we obtain a section of 
$\Lambda^2E\otimes K^4$. But from (\ref{decomp}) 
$$\varphi=\frac{1}{\sqrt{(v_0,v_0)}}i_{v_0}\rho.$$
Now $\varphi(u_1,u_2)=\omega(\Phi'Iu_1,u_2)$, and $I=\KK_{\Omega}/\sqrt{-\lambda(\Omega)}$.  Since  $\lambda(\Omega) = (v_0,v_0)$ it follows that, restricted to $V_1^*$, 
$i_{v_0}\rho(u_1,u_2)=\omega(\Phi'\KK_{\Omega}u_1,u_2).$
Now on $D$, $v_0$ is null and so lies in $V_1^*$, hence 
\begin{equation}
i_{v_0}\Omega(u_1,u_2) =\omega(\Phi'\KK_{\Omega}u_1,u_2) 
\label{Domega}
\end{equation}
 On $D$ we have  the normal form (\ref{norm2}) 
 $$\Omega=\xi_1\wedge \eta_2\wedge \eta_3+\xi_2\wedge \eta_3\wedge \eta_1+\xi_3\wedge \eta_1\wedge \eta_2$$
 where it is clear that $i_v\Omega=0$ if and only if $v=0$, so  Equation \ref{Domega} uniquely determines $v_0$ on $D$. Since this inclusion  is what we used to construct the bundle $V$ from $E$ as an extension in Section \ref{duality}, it is clear that in the $G_2$ case we do not have to consider a covering of the Prym variety as in the general $SO(2m+1)$ Higgs bundle. What we should check, however is that, starting from the symplectic bundle which defines the right hand side of Equation \ref{Domega}, there is a $v_0$ which satisfies the equation.
\vskip .25cm
We start then with $E$ and  $\Phi'$ preserving the symplectic form and $\Omega$. At a point on $D$ it is given by $\phi_0$, which lies in the Lie algebra of $G$.

 Now from the normal form (\ref{fiform}) of $\Omega$  we find that  $\KK_{\Omega}$ is given by $k_0$ where (with a standard trivialization of $\Lambda^6 E$)
$k_0(x,y)=(-2y,0)$. Thus, from (\ref{fiform}),  
$$\omega(\phi_0 k_0(x_1,y_1),(x_2,y_2))=-2(a\times y_1,y_2)=2(a,y_1\times y_2)$$
But 
$$\Omega((a,0),(x_1,y_1),(x_2,y_2))= (a,y_1\times y_2)$$
so $v_0=2a$ solves the equation.
\vskip .25cm
We then have the following
\begin{thm} Let $S$ be a curve of the form (\ref{x}) and $L_0$ be a line bundle in $P(S,C)^-$. Let $(E,\Phi')$ be the corresponding $Sp(6)$ Higgs bundle. Then the canonical vector $v_0$ in (\ref{Domega}) defines, as in Section \ref{reconstruct}, an extension $V$ which is a Higgs bundle with $G_2$ structure.
\end{thm}

\begin{prf} The line bundle $L_0$ is in $P(S,\bar S)$ and so defines a symplectic bundle $E=p_*U$. We define $\Omega\in H^0(\Sigma,\Lambda^3E^*\otimes K^{3/2})$ by push-down:  if $U_{\alpha}\subset \Sigma$ is an open set, sections of $E$ over $U_{\alpha}$ are sections of $U$ over $p^{-1}(U_{\alpha})$ which are sections of $p_{S\,*}U$ on $C$, and this bundle, because $L_0$ is in the Prym variety $P(S,C)$, has a twisted volume form which we evaluate on the three sections. We then obtain a Higgs bundle $(E,\Phi')$ where $\Omega$ is $\Phi'$-invariant. 

What remains is to show that the rank $7$ bundle obtained from the canonical extension admits a three-form $\rho$ which is everywhere in the open orbit. We adopt the point of view of Theorem \ref{canon} and in the orthogonal decomposition $(v,s)\in V_1^*\oplus K^{-m}$ outside of $D$ use the expression (\ref{decomp}) for $\rho(v_1+s_1, v_2+s_2,v_3+s_3)$. This gives
$$\Omega(v_1,v_2,v_3)+
\omega(\Phi'\KK_{\Omega}v_1,v_2)s_3+\omega(\Phi'\KK_{\Omega}v_1,v_2)s_3+\omega(\Phi'\KK_{\Omega}v_1,v_2)s_3.$$
Now, as before, write this relative to a local splitting where 
$$v= w - \frac{t}{cz}(v_0+zv_{01}+\dots),\qquad s=\frac{t}{cz}+\dots$$
and $v_0$ is the canonical vector. Evaluating this on vectors of this form will be smooth so long as 
$$\Omega(v_0,w_2,w_3)-\omega(\phi_0 k_0 w_2,w_3)=0$$
for all $w_2,w_3$. 
But this is the relation (\ref{Domega}).

It follows that we have $\rho\in H^0(\Sigma,\Lambda^3V^*)$ which extends our definition outside $D$. Now since $\Lambda^7V^*$ is trivial $\lambda(\rho)$ is a constant. It is non-zero since by construction it was non-zero outside $D$. At each point of $\Sigma$ it therefore lies in the open orbit and defines a $G_2$ structure on $V$. 

The $SO(7)$ Higgs field constructed in Section \ref{reconstruct} annihilated $v_0 \in H^0(\Sigma, V\otimes K^3)$. Since $\Phi'$ preserved $\Omega$ and $\omega$, $\Phi$ clearly preserves $\rho$ which is constructed out of these and we have a $G_2$ Higgs bundle. 

\end{prf}

\subsection{The abelian variety}

We have seen how a line bundle in the subgroup $P(S,C)^-\subset P(S,C)$ defines a $G_2$ Higgs bundle. To discuss duality we need to know more about this, and in particular that it is connected.  

First, let us calculate its dimension. If $TP$ is the tangent space to $P(S,C)$ at the origin then 
$H^1(S,{\mathcal O})\cong p^*_S H^1(C,{\mathcal O})\oplus TP.$ 
The involution $\sigma$ commutes with $p_S:S\rightarrow C$, so the anti-invariant parts satisfy 
$$H^1(S,{\mathcal O})^-\cong p_S^*H^1(C,{\mathcal O})^-\oplus TP^-.$$ 
This gives
\begin{eqnarray*}
\dim P(S,C)^-&=&(g(S)-g(\bar S))-(g(C)-g(\Sigma))\\
&=& (36(g-1)+1)-(15(g-1)+1)-(8(g-1)+1-g)\\
&=& 14(g-1)
\end{eqnarray*}
and this is $\dim G_2 (g-1)$ as expected.

\begin{prp} \label{connect}  $P(S,C)^-$ is connected.
\end{prp}
\begin{prf}
Note the names of the various projections:
$$\pi:S\rightarrow \bar S\qquad p_C: C\rightarrow \Sigma\qquad p_S:S\rightarrow C \qquad \pi_{\bar S}:\bar S\rightarrow \Sigma.$$
We write the group law additively here. Let $A$ be the identity component of $P(S,C)^-$. Since $P(S,C)$ is connected   $x\mapsto x-\sigma x$ maps $P(S,C)$ onto $A$. For $x\in P(S,C)^-$ take $y\in P(S,C)$ such that $x=2y$ and write
\begin{equation}
x=y+\sigma y+y-\sigma y
\label{xyz}
\end{equation}
Then $z=y+\sigma y$ is pulled back from $\bar S$ and  satisfies $\sigma z=-z$ so $z=-z$ and  lies in $\pi^*H^1(\bar S,\Z_2)$.
\vskip .25cm
 Consider the endomorphism $s$  defined by $s(x)=2x$ on $A$. We have seen that there is a canonical choice of extension to define $V$, so this means, comparing with the $SO(2m+1)$ case in Section \ref{solie}, that there is a section of 
$s: A/(A\cap \pi^*H^1(\bar S,\Z_2))\rightarrow A$
or equivalently, $$A_2= \pi^*H^1(\bar S,\Z_2)\cap A$$ (where the subscript $2$ denotes the elements of order $2$).

The map $\pi_{\bar S}:\bar S\rightarrow \Sigma$ is of degree $3$ so given $y\in H^1(\Sigma,\Z_2)$ we can write $y=3y=\Nm_{\pi_{\bar S}}\pi^*_{\bar S}y$ for an element of order $2$, and this gives a decomposition $x\mapsto (x+\pi^*_{\bar S}\Nm_{\pi_{\bar S}}x, \Nm_{\pi_{\bar S}}x)$ 
$$H^1(\bar S,\Z_2)\cong P(\bar S,\Sigma)_2\oplus \pi_{\bar S}^*H^1(\Sigma,\Z_2).$$
Now the order of $P(\bar S,\Sigma)_2$ is
$$2^{2(15(g-1)+1-g)}=2^{2(14(g-1))}$$
which is the order of $A_2$. Moreover if $u\in H^1(\bar S,\Z_2)$ and $\pi^*u\in P(S,C)$ then 
$\Nm_{p_S}\pi^*u=0$. But $S$ is the fibre product of $\pi_{\bar S}:\bar S\rightarrow \Sigma$ and $p_C:C\rightarrow \Sigma$ hence 
$$0=\Nm_{p_S}\pi^*u=p_C^*\Nm_{\pi_{\bar S}}u.$$
Since $p_C^*$ is injective $\Nm_{\pi_{\bar S}}u=0$ and so $u\in P(\bar S,\Sigma)_2$. We deduce that $A_2=\pi^*P(\bar S,\Sigma)_2$.
\vskip .25cm
Now $y+\sigma y$ in (\ref{xyz}) is of order $2$ and of the form $\pi^*u$ and lies in $P(S,C)$. It follows that $y+\sigma y\in A_2$. Thus $x=y+\sigma y+y-\sigma y\in A$ is a sum of two elements of $A$ and so  $P(S,C)^-=A$.
\end{prf}

\section{Duality for $G_2$}
\subsection{The dual variety}

\begin{prp} \label{dualp} The dual of the abelian variety $P(S,C)^-$ is 
$$P(S,C)^-/p_S^*H^1(C,\Z_3)^-.$$
\end{prp}
\begin{prf} The abelian variety $P(S,C)^-$  is the kernel of 
$\Nm_{\pi}$ restricted to $P(S,C)$, and since $\Nm_{\pi_{\bar S}}\Nm_{\pi}=\Nm_{p_C}\Nm_{p_S}$, its image is contained in  $P(\bar S, \Sigma)$. The dual of $P(S,C)$ is $P(S,C)/p_S^*H^1(C,\Z_3)$ and there is a surjective homomorphism from this group to  $(P(S,C)^-)^{\vee}$, since $P(S,C)^-\subset P(S,C)$ is connected. Restricting to the anti-invariant part gives a surjection
$$P(S,C)^-/H^1(C,\Z_3)^-\rightarrow (P(S,C)^-)^{\vee}.$$
The kernel of this is the image of the dual of $P(\bar S, \Sigma)$, which is 
$P(\bar S, \Sigma)/\pi_{\bar S}^*H^1(\Sigma,\Z_3)$. 

But $P(S,\bar S)$ intersects $\pi^*J(\bar S)$ in elements of order $2$ and $P(S,C)^-=P(S,\bar S)\cap P(S,C)$. Hence $\pi^*P(\bar S, \Sigma)\cap P(S,C)^-\subset \pi^*P(\bar S, \Sigma)_2.$ But we saw in the proof of Proposition \ref{connect} that this consists of {\it all} elements of order $2$ in $P(S,C)^-$.  Because $2$ and $3$ are coprime, it follows that the dual $(P(S,C)^-)^{\vee}$ is the quotient of $P(S,C)^-/\pi^*H^1(C,\Z_3)^-$ by all elements of order $2$ and $x\mapsto 2x$ identifies this with itself.
\end{prf}

We shall find this variety appearing as the abelian variety for a different fibre in the Higgs bundle moduli space. 

\subsection{The cameral curve}\label{camera}

The spectral curve $S$ is a $6$-fold cover of $\Sigma$. Its equation is a cubic in $x^2$ whose discriminant is 
\begin{equation}
\Delta=q(\frac{1}{2}f^3-27q)=27q q^{\vee}.
\label{discr}
\end{equation}
where
$$q=(\lambda_1\lambda_2\lambda_3)^2,\quad 27q^{\vee}=((\lambda_1-\lambda_2)(\lambda_2-\lambda_3)(\lambda_3-\lambda_1))^2.$$
Now by definition, $S$ is a curve on which $x$ is a single valued eigenvalue of $\Phi'$. Thus on $S$ we can find the other eigenvalues by fully factorizing the polynomial
$$(w-x^2)(w^2+bw+c)=w^3-fw^2+\frac{f^2}{4}w-q.$$
Here $b=x^2-f$ and $c=(x^2-f/2)^2$ and we calculate the discriminant $b^2-4c$ of the quadratic factor to be
$x^2(2f-3x^2)$. So we can solve the quadratic by setting
\begin{equation}
3y^2=2f-3x^2
\label{y}
\end{equation}
to obtain 
$$w=\frac{1}{2}(-b\pm \sqrt{3}xy)=-\frac{1}{4}(x^2+3y^2\mp 2\sqrt{3}xy).$$
The six roots $\pm \lambda_i$ of the equation are therefore
\begin{equation}
\lambda_1= x,\quad \lambda_2 = (-x+\sqrt{3}y)/2,\quad \lambda_3 = (-x-\sqrt{3}y)/2.
\label{eqeigen}
\end{equation}
\vskip .25cm
To get to this point, we introduced the curve $W$ given by $3y^2=2f-3x^2$. 
It lies in the three-dimensional manifold $K\otimes \C^2\rightarrow \Sigma$ and is given by the two equations (\ref{x}), (\ref{y})
\begin{equation}
x^2+y^2=2f/3,\quad x^6-fx^4+\frac{f^2}{4}x^2-q=0.
\label{Sequation}
\end{equation}
It is a double covering of $S$ branched over $y=0$, a divisor of $p^*K$. Since $K_S\cong p^*K^6$  this means that $K_W$ is the pullback of $K^7$. Hence
$$2g(W)-2=12\times 7 \times 2g-2$$
and $g(W)=84(g-1)+1.$

There is an action of the dihedral group $D_6$ of order $12$ on $W$: firstly a rotation $r$ by $\pi/3$ is given by the matrix
$$\pmatrix {1/2 &\sqrt{3}/2\cr
                                -\sqrt{3}/2& 1/2}$$
                                and this maps
                                $$\pmatrix{x\cr y}\mapsto \pmatrix{(x+\sqrt{3}y)/2\cr(y-\sqrt{3}x)/2 }\mapsto \pmatrix{(-x+\sqrt{3}y)/2\cr(-y-\sqrt{3}x)/2 }\mapsto \pmatrix{-x\cr -y}$$
so the first entry runs through the six eigenvalues in (\ref{eqeigen}). Together with the reflection $s$ defined by $(x,y)\mapsto (x,-y)$ which  defines the double covering $W\rightarrow S$, this generates the $D_{6}$ action:  
$s^2=1,r^6=1$ and $rs=sr^{-1}.$

\vskip .25cm
Substituting for $x^2$ in (\ref{Sequation}) gives the equivalent formulation:
\begin{equation}
x^2+y^2=2f/3,\quad y^6-fy^4+\frac{f^2}{4}y^2+q-\frac{f^3}{54}=0
\label{Svequation}
\end{equation}
so that replacing $q$ by $q^{\vee}$ gives a different spectral curve $S^{\vee}$ with the same  curve $W$. Duality for $G_2$ entails interchanging the roles of $S$ and $S^{\vee}$.

\begin{rmk} The dihedral group $D_{6}$ is the Weyl group of $G_2$ and $W$ is then the {\it cameral curve} of $\Sigma$ discussed in the root system treatment in \cite{Don}.
\end{rmk}

\subsection{Dual curves} 

The spectral curve $S$ is the quotient of $W$ by the reflection  $s(x,y)=(x,-y)$, and $S^{\vee}$ the quotient by $r^3s(x,y)=(-x,y)$. There are two conjugacy classes of reflections in this dihedral group --  reflections in an axis passing through two opposite vertices of a hexagon, and those in an axis through the mid-points of opposite sides. The reflection $s$ belongs to one and $r^3s$ to the other. But $rs$ is conjugate to $r^3s$, so the curve $S^o$ defined as the quotient of $W$ by $rs$, is isomorphic to $S^{\vee}$. The intermediate curve $C$ is the quotient of $W$ by the $D_3$ generated by $r^2, s$, and there is a corresponding curve $C^o$ for the group generated by $r^2, rs$.  We shall relate the abelian variety for $S^o$ to the dual of  the variety for $S$.

Let $f:W\rightarrow S$ and $f^0:W\rightarrow S^o$ be the quotient maps, then:
\begin{prp}
$\Nm_{f^o}f^*$ defines an isomorphism from $P(S,C)^-/p_S^*H^1(C,\Z_3)^-$ to $P(S^o,C^o)^-$.
\end{prp}
Using Proposition \ref{dualp}, this result shows that the dual of $P(S,C)^-$ is isomorphic to $P(S^{o},C^{o})^-$. Together with $P(S^o,C^o)^-\cong P(S^{\vee}, C^{\vee})^-$ this realizes Langlands duality within the same moduli space. Note that if we simply  pull back from $S$ and push down  to $S^{\vee}$ we get zero, which is why we use $S^o$ instead of $S^{\vee}$.

\begin{prf} The result follows from the more general results  of Carocca et al \cite{Car}. We tailor their method here to our specific situation. 

First consider  the quotient of $W$ by the subgroup $\Z_3$  generated by  $r^2$ and denote by  $\pi_W:W\rightarrow W/\Z_3$ the quotient map. The curve $W/\Z_3$ is a ramified  double cover of $C=W/D_3$ with projection $g$. The spectral curve  $S$ is the quotient of $W$ by the reflection $s\in D_3$ with projection $f:W\rightarrow S$. Then the curve $W$ may be considered as the fibre product of 
$g:W/\Z_3\rightarrow C$ and $p_S:S\rightarrow C$. In particular
$$\Nm_{\pi_W} f^*=g^*\Nm_{p_S}.$$ 

If $x\in P(S,C)$, then $\Nm_{p_S}x=0$ and 
$$0=g^*\Nm_{p_S}x=\Nm_{\pi_W} f^*x.$$
Let $y=f^*x$, then this means that  $(1+r^2+r^4)y=0$ and $sy=y$. 
\vskip .25cm
 Now  suppose $x$ lies in the kernel of $\Nm_{f^o}f^*$. Then $(1+rs)y=0$. But $sy=y$ so $ry=-y$ and $(1+r^2+r^4)y=0$ gives $3y=0$.
\vskip .25cm
Now we have $r^2y=y$ and $sy=y$ so  $y$ is invariant under the dihedral group $D_3$. It is the class of a line bundle pulled back from $C=W/D_3$ if the action  at the fixed points of elements in the group is trivial. Now a  rotation in $D_{6}$ only fixes the origin in $K\otimes \C^2$ and this is $x=y=0$.  In the generic case,  this does not lie on the curve $W$, so there are no fixed points for $r^2$. On the other hand, $y=f^*x$ and so the action at fixed points of $s$ is trivial. But in $D_3$, any two reflections are conjugate, so the action  is trivial at all fixed points and therefore $y$ is pulled back from $C$, and hence 
 $x\in p_S^*H^1(C,\Z_3).$
 
If $\sigma x=-x$ then $x\in  p_S^*H^1(C,\Z_3)^-$. Conversely if $x\in  H^1(C,\Z_3)^-$, then $y=f^*x$ is invariant under $D_3$ so $sy=y, r^2y=y$. Since $\sigma x=-x$, $r^3y=-y$ and so $ry=-y$. This means that $(1+rs)y=(1+r)y=0$ and $x$ is in the kernel of $\Nm_{f^o}f^*$.
 
\end{prf}
\begin{rmk} We have seen that  the involution
$$(f,q)\mapsto(f, \frac{f^3}{54}-q)$$ 
on $H^0(\Sigma, K^2)\oplus H^0(\Sigma, K^6)$ takes a fibre to its dual. When $f=0$, the two spectral curves $x^6\pm q=0$ are isomorphic and one might expect the abelian variety to be dual to itself.  This is indeed the case: $S$ has an action of $\Z_6$ generated by $rx= e^{i\pi/3}x$ and $C$ is the quotient by the $\Z_3$ generated by  $r^2$. The map $z\mapsto (1+r)z$ of $P(S,C)^-$ to itself has kernel $p_S^*H^1(C,\Z_3)^-$.
\end{rmk}

\subsection{The $D_6$ action}

The pull back $f^*P(S,C)$ to the curve $W$ is characterized by the condition $sx=x$ and $(1+r^2+r^4)+s(1+r^2+r^4)x=0$ since $C$ is the quotient of $W$ by the group $D_3$ generated by $s$ and $r^2$. The anti-invariant part $P(S,C)^-$ satisfies the further condition $r^3x=-x$. Thus its  tangent space $T\subset H^1(W,{\mathcal O})$ is the solution to the equations  
\begin{equation}
(1+r^2+r^4)x=0\,\,\,\,\,\, r^3x=-x
\label{rr}
\end{equation}
and $sx=x$. Similarly $T^{\vee}$, the tangent space of $(f^v)^* P(S^{\vee},C^{\vee})$ satisfies (\ref{rr}) and $sx=-x$.

Since $sr=r^{-1}s=-r^2s$, if $x\in T$ then
$$rx=\frac{1}{2}(rx+srx)+\frac{1}{2}(rx-srx)=\frac{1}{2}(rx-r^2x)+\frac{1}{2}(rx+r^2x)$$
and both factors satisfy (\ref{rr}) so that $rx$ lies in $T\oplus T^{\vee}$. It follows that the $28(g-1)$-dimensional space $T\oplus T^{\vee}$ is preserved by the $D_6$ action. Moreover the relations above show that this is the subspace of $H^1(W,{\mathcal O})$ whose  isotype  is the  two-dimensional irreducible dihedral representation. Equivalently 
$$T\oplus T^{\vee}=\C^2\otimes V$$
for some $14(g-1)$-dimensional vector space $V$.

 The pull-back of $P(S,C)^-$ and $P(S^{\vee},C^{\vee})^-$ generate a $28(g-1)$-dimensional abelian variety in $H^1(W,{\mathcal O}^*)$ on which $D_6$ acts. Although their tangent spaces $T$ and $T^{\vee}$ are complementary, the abelian variety is not a product, because there is a non-zero intersection. In fact if $x\in  f^*P(S,C)^-\cap (f^{\vee})^*P(S^{\vee},C^{\vee})^-$ then $sx=x=-x$ and $x$ is of order $2$. But in the proof of Proposition \ref{connect} we saw that the group of elements of order $2$ in $P(S,C)^-$ is $\pi^*P(\bar S,\Sigma)_2$.  Here $\bar S$ is the quotient of  $S$ by the involution, which is the quotient of $W$ by the group $1, r^3, s, r^3s$. But $S^{\vee}$ is the quotient of $W$ by  $r^3s$ so  $\bar S=\overline{ (S^{\vee})}$. There is thus a natural identification of the elements of order $2$ in $P(S,C)^-$ and  $P(S^{\vee},C^{\vee})^-$ and the abelian variety is the quotient by the diagonal action. The squaring map on either factor defines a homomorphism to $P(S,C)^-$ with kernel $P(S^{\vee},C^{\vee})^-$ or vice-versa.

\begin{rmk} Donagi's root system approach to spectral curves describes the abelian variety as the identity component of the moduli  space of  Weyl-invariant  $H$-bundles on the cameral curve, where $H$ is the Cartan subgroup.  For $G_2$, the Cartan subalgebra is  $\C^2$ with the Weyl group action the dihedral representation. As we have seen, the $D_6$-invariant part of  $\C^2\otimes H^1(W,{\mathcal O})$ is isomorphic to $V$, and this can be identified with the $s=1$ subspace of $\C^2\otimes V$ which is $T=TP(S,C)^-$. It follows from this and the connectedness that  our description of the $G_2$ abelian variety and that of Donagi coincide.
\end{rmk}

\subsection{The cubic form}
An open set (the complement of the discriminant locus) of the base space $B$ of an algebraically completely integrable  Hamiltonian system  has a natural differential geometric structure on it called a {\it special K\"ahler structure} (see \cite{Freed}, \cite{H5}). This involves distinguished flat coordinates ({\it not} the flat vector space coordinates for our integrable system) and a cubic form -- a holomorphic section of $Sym^3 T^*$ (introduced initially by Donagi and Markman \cite{Don}). In fact in flat coordinates, the cubic form is the third derivative of a holomorphic function.  

 For our $G_2$ Higgs bundle moduli space, we have an involution 
 $$(f,q)\mapsto (f,q^{\vee})=(f,\frac{1}{54}f^3-q)$$
 on  $B=H^0(\Sigma, K^2)\oplus H^0(\Sigma, K^6)$ and it seems quite likely that this is an isometry of the special K\"ahler structure.  We shall restrict ourselves here to calculating the cubic form, using recent work of Balduzzi \cite{Bald} and show that this is invariant under the involution. 

The cubic form is essentially the infinitesimal period map. A tangent vector  $u\in T_b$, the tangent space of $B$  at $b$, defines a   Kodaira-Spencer deformation class in $H^1(X_b,T)$ where $X_b$ is   the fibre over $b$.  The cup product  gives  a linear map $\chi_u: H^0(X_b,T^*)\rightarrow H^1(X_b,{\mathcal O})$, or $\chi_u\in Sym^2 H^0(X_b,T^*)^*$. The symplectic form on the total space  identifies $H^0(X_b,T^*)$ with $T_b$, and then $\chi_u(v, w)$ is the cubic form. 

Building on unpublished work of Pantev, Balduzzi has given a formula for the cubic form where the integrable system is the Higgs bundle moduli space. He identifies the tangent space at $b$ in the base as the space of Weyl-invariant sections of $\lie{h}\otimes K_W$ on the cameral curve $W$. The formula is \cite{Bald}

\begin{equation}
\chi_u(v,w)=\sum_{D(a)=0}\Res_a^2 \frac{D_u}{D} B(v,w).
\label{cube}
\end{equation}

 Here $B$ is the Killing form and $B(v,w)$ is a quadratic differential on the cameral curve. The discriminant locus on $\Sigma$ is given by $D$, a section of $K^n$ where $n$ is the order of the Weyl group. This section is a polynomial in the differentials $\bigoplus_1^n H^0(\Sigma,K^{d_i})$ which form the base of the fibration, and $D_u/D$ is the logarithmic derivative in the direction $u$. The expression $\Res^2_a(q)$ of a quadratic differential is the coefficient of $dw^2/w^2$ in a local coordinate with $w(a)=0$. The definition is written  in terms of the cameral curve but the final expression is well-defined on $\Sigma$.
 \vskip .25cm
 We  described in the previous section the Weyl-invariant elements  in $\lie{h}\otimes H^1(W,{\mathcal O})$ for $G_2$. We now want the invariant subspace of $\lie{h}\otimes H^0(W,K_W)$. This is naturally isomorphic to the tangent space of $T_b$ of the base, which is the space of infinitesimal deformations of the spectral curve  
$$x^6-fx^4+\frac{f^2}{4}x^2-q=0.$$
Let $(\dot f,\dot q)\in H^0(\Sigma, K^2)\oplus H^0(\Sigma, K^6)$ denote such a deformation, then we consider the section 
$$-\dot f x^4+\frac{f \dot f}{2}x^2-\dot q$$
of $p^*K^6$ on $S$ which is the first order deformation of the equation.
This is pulls back to a differential on $W$, and we saw in Section \ref{camera} that $K_W\cong (pf)^*K^7$. Using this isomorphism, the  differential on $W$ is 
$$X=-\dot f x^4y+\frac{f \dot f}{2}x^2y-\dot q y.$$
 We can do the same for the other spectral curve $S^{\vee}$ 
 $$y^6-fy^4+\frac{f^2}{4}y^2 -q^{\vee}=y^6-fy^4+\frac{f^2}{4}y^2+q-\frac{f^3}{54}=0$$
 and get a differential
 $$Y=-\dot f y^4x+\frac{f \dot f}{2}y^2x+\dot q x-\frac{f^2\dot f}{18}x.$$
 We claim that the space of such pairs
 $$(X,Y)  = \dot f(- x^4y+\frac{f}{2}x^2y, - y^4x+\frac{f }{2}y^2x-\frac{f^2}{18}x)+\dot q\left(-y,x\right)$$
transforms according to the $D_6$ dihedral representation and thus consists  of the Weyl-invariant ${\lie h}$-valued differentials. 

 This is easier to see by using the relation $x^2+y^2=2f/3$ and putting $x=\sqrt{2f/3}\cos \theta$ and $y= \sqrt{2f/3} \sin \theta$, for then the above expression simplifies to 
$$\dot f\frac{1}{36}\sqrt{\frac{2f}{3}}f^2(-\sin 5\theta+\sin \theta, - \cos 5\theta -\cos \theta)+\dot q\sqrt{\frac{2f}{3}}(-\sin \theta,\cos \theta).$$
Now use 
$$\dot q^{\vee}=\frac{\dot f f^2}{18}-\dot q$$ to write this as 
$$\sqrt{{f}/{6}}[\dot q(-\sin \theta-\sin 5\theta,\cos \theta-\cos 5\theta)+\dot q^{\vee}(\sin \theta-\sin 5\theta, -\cos \theta-\cos 5\theta)].$$
Applying the inner product $B$, which is just  the Euclidean inner product on $\C^2$, we get the quadratic expression $B(v,w)$ in formula (\ref{cube})
$$\frac{f}{3} [(1-\cos 6\theta)\dot q_1 \dot q_2+ (1+\cos 6\theta) \dot q_1^{\vee} \dot q_2^{\vee}]$$ where $v=(\dot f_1,\dot q_1), w=(\dot f_2,\dot q_2)$. 

Now use $\cos 6\theta=32 \cos^6\theta-48 \cos ^4\theta+18 \cos^2\theta-1$ and the equation of the spectral curve, and we obtain 
\begin{equation}
B(v,w)=36\frac{1}{f^2}(q^{\vee}\dot q_1 \dot q_2+q\dot q_1^{\vee} \dot q_2^{\vee}).
\label{Bvw}
\end{equation}
From (\ref{discr}) the discriminant divisor is given by the section $q q^{\vee}$ of $K^{12}$. Generically $q$ and $q^{\vee}$ have disjoint zeros  so the cubic form (\ref{cube}) is in this case the sum of two terms 
$$36\sum_{q(a)=0}\Res_a^2 \frac{q^{\vee}}{qf^2}\dot q_1\dot q_2\dot q_3+36\sum_{q^{\vee}(a)=0}\Res_a^2 \frac{q}{q^{\vee}f^2}\dot q^{\vee}_1\dot q^{\vee}_2\dot q^{\vee}_3$$
which is clearly invariant under the involution $(f,q)\mapsto (f,q^{\vee})$.

We need to write this in terms of a local coordinate on $W$ to evaluate the residues at the zeroes of $q$ and $q^{\vee}$. In fact, since $W$ is the double covering of $S$ branched over $q^{\vee}=0$, and $q$ has no common zeroes with $q^{\vee}$, we can evaluate at the zeroes of $q$ using a coordinate on the spectral curve $S$. 

Let $z$ be a local coordinate on $\Sigma$, so that $f=g(z)dz^2$ and $q=r(z)dz^6$. The tautological section $x$ of $p^*K$ on $K$ is then just $wdz$, and the spectral curve  has equation  
$$w^6-g(z)w^4+\frac{g(z)^2}{4}w^2-r(z)=0$$
and since $r'(a)$ is nonzero where $r(a)=0$, $w$ is a local coordinate on $S$ near $a$.  

  At $w=0$, $r'(z)dz$ is a nonvanishing section of $p^*K^6\otimes N^*$, where $N$ is the normal bundle  of $S\subset K$. The canonical one-form on the cotangent bundle of $\Sigma$ is $wdz$ and its derivative $dw\wedge dz$ is the symplectic form on $K$. We use this to identify the canonical bundle $K_S$  with $p^*K^6$, so a section $s$ of $p^*K^6$ defines a differential with the local form $sdw/r'$ on $S$. On $W$, where $K_W$ is the pullback of $K^7$, this corresponds to  
$$sy\frac{dw}{r'}.$$

Where $q$ vanishes, only the first term in (\ref{Bvw})  contributes to the residue, and where $q=0$, $q^{\vee}=f^3/54$. So this term is 
$$36\frac{1}{f^2}q^{\vee}\dot q_1 \dot q_2=\frac{2}{3}f\dot q_1 \dot q_2=y^2\dot q_1 \dot q_2$$
since $y^2=2f/3$ where $x=0$. This quadratic differential thus has the local form
$$\dot r_1(z) \dot r_2(z)\frac{dw^2}{r'(z)^2}.$$

Multiplying by $\dot q_3/q$ and using $q=f^2x^2/4+\dots$ gives 
$$\frac{\dot r_1\dot r_2\dot r_3}{g^2r'^2}(a)\frac{dw^2}{w^2}+\dots$$
which determines the residue term.

We can write this invariantly on $\Sigma$ now, since at a zero $a$ of the section $q$ of $K^6$, $q'=r'dz^7$ is a well-defined vector in $K^7_a$. Taking  $\dot q_i$  in $K^6_a$ and $f=gdz^2$  in $K^2_a$, we obtain  
$$\frac{\dot q_1\dot q_2\dot q_3}{f^2q'^2}(a)$$
which is simply a complex number.

Taking into account the double covering $W\rightarrow S$, the final formula for the cubic form is  
$$2\sum_{q(a)=0}\frac{\dot q_1\dot q_2\dot q_3}{f^2q'^2}(a)+2\sum_{q^{\vee}(a)=0}\frac{\dot q^{\vee}_1\dot q^{\vee}_2\dot q^{\vee}_3}{f^2(q^{\vee})'^2}(a)$$

\vskip 1cm
 Mathematical Institute, 24-29 St Giles, Oxford OX1 3LB, UK
 
 hitchin@maths.ox.ac.uk
\end{document}